%% file: apostextlaplacenum3.tex
\documentclass[12pt,a4paper,twoside]{article}

\title{\sc Conforming and Non-Conforming\\
Functional A Posteriori Error Estimates\\
for Elliptic Boundary Value Problems\\
in Exterior Domains:\\
Theory and Numerical Tests}
\def\shorttitle{A Posteriori Error Estimates for Elliptic Problems in Exterior Domains}
\def\pauthor{Olli Mali, Alexey Muzalevskiy and Dirk Pauly}

\def\mylabelonoff{off}
\def\allowdisbrk{no}

\input head_paper

\input head_paule_olli

\begin{document}

\maketitle{}

\begin{center}
{\tt Dedicated to Sergey Igorevich Repin\\
on the occasion of his 60th birthday}
\end{center}

\begin{abstract}
This paper is concerned with the
derivation of conforming and non-conforming
functional a posteriori error estimates
for elliptic boundary value problems in exterior domains.
These estimates provide computable and guaranteed upper and lower bounds
for the difference between the exact and the approximate solution
of the respective problem. We extend the results from
\cite{paulyrepinell} to non-conforming approximations,
which might not belong to the energy space and are just considered
to be square integrable. Moreover, we present some numerical tests.\\
\end{abstract}

\tableofcontents

\section{Introduction}

As in \cite{paulyrepinell},
we consider the standard elliptic Dirichlet boundary value problem
\begin{align}
-\div A\nabla u&=f&&\text{in }\om,\label{diffeq}\\
u&=u_0&&\text{on }\ga,\label{bceq}
\end{align}
where $\om\subset\rN$ with $N\geq3$
is an exterior domain, i.e., a domain with compact complement,
having for simplicity a Lipschitz continuous boundary $\ga:=\pom$.
Moreover, $A:\om\to\rNtN$ is a real and symmetric $\Liom$-matrix valued function such that
$$\exists\,\alpha>0\quad\forall\,\xi\in\rN\quad\forall\,x\in\om\qquad
A(x)\xi\cdot\xi\geq\alpha^{-2}|\xi|^2$$
holds. As usual, when working with exterior domain problems,
we use the polynomially weighted Lebesgue spaces
$$\Ltsom:=\set{u}{\rho^{s}u\in\Ltom},\quad\rho:=(1+r^2)^{1/2}\cong r,\quad s\in\rz,$$
where $r(x):=|x|$ is the absolute value. Throughout the paper at hand
we just need the values $s\in\{-1,0,1\}$ of weights.
If $s=0$, then we write $\Ltom:=\Ltzom$.
Moreover, we introduce the polynomially weighted Sobolev spaces
\begin{align*}
\Homoom&:=\set{u}{u\in\Ltmoom,\,\na u\in\Ltom},\\
\Dom&:=\set{v}{v\in\Ltom,\,\div v\in\Ltpoom},
\end{align*}
which we equip as $\Ltsom$ with the respective scalar products.
We will not distinguish in our notation between scalar and vector valued spaces.
Moreover, to model homogeneous boundary traces we define as closure of test functions
$$\Hocmoom:=\ol{\Cicom}^{\Homoom}.$$
Note that all these spaces are Hilbert spaces and we have for the norms
\begin{align*}
\norm{u}_{\Ltsom}^2&=\norm{\rho^su}_{\Ltom}^2=\int_{\om}(1+r^2)^{s}|u|^2\dl,\\
\norm{u}_{\Homoom}^2&=\norm{\rho^{-1}u}_{\Ltom}^2+\norm{\na u}_{\Ltom}^2,\\
\norm{v}_{\Dom}^2&=\norm{v}_{\Ltom}^2+\norm{\rho\div v}_{\Ltom}^2.
\end{align*}
Also, let us introduce
for vector fields $v\in\Ltom$ the weighted norm
\begin{align}
\mylabel{defscp}
\norm{v}_{\Ltom,A}:=\scp{v}{v}_{\Ltom,A}^{1/2}
:=\scp{Av}{v}_{\Ltom}^{1/2}=\norm{A^{1/2}v}_{\Ltom}.
\end{align}
Let
$$\cN:=\frac{2}{N-2},\quad\cNa:=\alpha\cN.$$
From \cite[p. 57]{leisbook} we cite the Poincar\'e estimate III
(see also the appendix of \cite{paulyrepinell})
\begin{align}
\mylabel{poincare}
\forall\,u\in\Hocmoom\qquad\norm{u}_{\Ltmoom}
\leq\cN\norm{\na u}_{\Ltom}
\leq\cNa\norm{\na u}_{\Ltom,A},
\end{align}
which is the proper coercivity estimate for the problem at hand.
Using this estimate it is not difficult to get 
by standard Lax-Milgram theory unique solutions
$\uh\in\Hocmoom+\{u_{0}\}$ of \eqref{diffeq}-\eqref{bceq}
depending continuously on the data
for any $f\in\Ltpoom$ and $u_{0}\in\Homoom$.
Note that the solution $\uh$ satisfies the variational formulation
\begin{align}
\mylabel{varform}
\forall\,u\in\Hocmoom\qquad
\scp{\na\uh}{\na u}_{\Ltom,A}=\scp{A\na\uh}{\na u}_{\Ltom}=\scp{f}{u}_{\Ltom},
\end{align}
where we use the $\Ltom$-inner product notation 
also for the $\Ltmoom$-$\Ltpoom$-duality. Note that 
$$\scp{f}{u}_{\Ltom}=\int_{\om}fu\dl$$ 
is well defined since the product $fu$ belongs to $\Loom$.
Moreover, we note that
$$A\na\uh\in\Dom,\quad-\div A\na\uh=f.$$

Let $\ut$ be an approximation of $\uh$. The aim of this contribution is twofold.
First, we extend the results from \cite{paulyrepinell} to non-conforming
approximations $\ut$ which no longer necessarily belong to the
natural energy space $\Homoom$ and hence lack regularity.
We will just assume that we have been given an approximation $\vt\in\Ltom$ of $A\na\uh$
without any regularity except of $\Ltom$.
Second, we validate the conforming a posteriori error estimates
for the problem \eqref{diffeq}-\eqref{bceq}
by numerical computations in the exterior domain $\om$
and therefore demonstrate that this technique also works in unbounded domains.
Such a posteriori error estimates
have been extensively derived and discussed earlier for problems in bounded domains, see e.g.
\cite{NeittaanmakiRepin2004,repinapostvarprob,repinbookone} and the literature cited there.
The underlying general idea is to construct estimates via Lagrangians. 
In linear problems this can be done by splitting
the residual functional into two natural parts
using simply integration by parts relations, 
which then immediately yield guaranteed and computable 
lower and upper bounds. In fact, one adds a zero to the weak form.

\section{Conforming A Posteriori Estimates}

For the convenience of the reader, we repeat also in the conforming case
the main arguments from \cite{paulyrepinell} to obtain the desired a posteriori estimates.
We want to deduce estimates for the error
$$e:=\uh-\ut$$
in the natural energy norm $\norm{\na e}_{\Ltom,A}$.
In this section we only consider conforming approximations, i.e.,
$\ut\in\Hocmoom+\{u_{0}\}$ and therefore $e\in\Hocmoom$.

Introducing an arbitrary vector field $v\in\Dom$ and inserting a zero into \eqref{varform},
we have for all $u\in\Hocmoom$
\begin{align}
\mylabel{calcone}
\begin{split}
\scp{\na(\uh-\ut)}{\na u}_{\Ltom,A}
&=\scp{f}{u}_{\Ltom}-\scp{A\na\ut-v+v}{\na u}_{\Ltom}\\
&=\scp{f+\div v}{u}_{\Ltom}-\scp{\na\ut-A^{-1}v}{\na u}_{\Ltom,A}
\end{split}
\end{align}
since $u$ satisfies the homogeneous boundary condition.
Therefore, by \eqref{poincare}
\begin{align}
\mylabel{calctwo}
\begin{split}
\scp{\na e}{\na u}_{\Ltom,A}
&\leq\norm{f+\div v}_{\Ltpoom}\norm{u}_{\Ltmoom}
+\norm{\na\ut-A^{-1}v}_{\Ltom,A}\norm{\na u}_{\Ltom,A}\\
&\leq\big(\ub{\cNa\norm{f+\div v}_{\Ltpoom}
+\norm{\na\ut-A^{-1}v}_{\Ltom,A}}_{\ds=:M_{+}(\na\ut,v;f,A)=:M_{+}(\na\ut,v)}\big)
\norm{\na u}_{\Ltom,A}.
\end{split}
\end{align}
Taking $u:=e\in\Hocmoom$ we obtain the upper bound
\begin{align}
\mylabel{confub}
\norm{\na e}_{\Ltom,A}&\leq\inf_{v\in\Dom}M_{+}(\na\ut,v).
\end{align}
We note that for $v:=A\na\uh$ we have $M_{+}(\na\ut,v)=\norm{\na e}_{\Ltom,A}$.
Therefore, in \eqref{confub} we even have equality.

The lower bound can be obtained as follows. 
Let $u\in\Hocmoom$. Then, by the trivial inequality
$\norm{\na(\uh-\ut) - \na u}_{\Ltom,A}^2 \geq 0$
and \eqref{varform} we obtain
\begin{align*}
\norm{\na(\uh-\ut)}_{\Ltom,A}^2
&\geq2\scp{\na(\uh-\ut)}{\na u}_{\Ltom,A}-\norm{\na u}_{\Ltom,A}^2\\
&=\ub{2\scp{f}{u}_{\Ltom}-\scp{\na(2\ut+u)}{\na u}_{\Ltom,A}}_{\ds=:M_{-}(\na\ut,u;f,A)=:M_{-}(\na\ut,u)}
\end{align*}
and hence we get the lower bound
\begin{align}
\mylabel{conflb}
\norm{\na e}_{\Ltom,A}^2&\geq\sup_{u\in\Hocmoom}M_{-}(\na\ut,u).
\end{align}
Again, we note that for $u:=\uh-\ut=e$ we have $M_{-}(\na\ut,u)=\norm{\na e}_{\Ltom,A}^2$.
Therefore, also in \eqref{conflb} equality holds.
Let us summarize:

\begin{theo}[conforming a posteriori estimates]
\mylabel{conftheo}
Let $\ut\in\Hocmoom+\{u_{0}\}$. Then
$$\max_{u\in\Hocmoom}M_{-}(\na\ut,u)
=\norm{\na(\uh-\ut)}_{\Ltom,A}^2
=\min_{v\in\Dom}M_{+}^2(\na\ut,v),$$
where the upper and lower bounds are given by
\begin{align*}
M_{+}(\na\ut,v)
&=\cNa\norm{f+\div v}_{\Ltpoom}
+\norm{\na\ut-A^{-1}v}_{\Ltom,A},\\
M_{-}(\na\ut,u)
&=2\scp{f}{u}_{\Ltom}
-\scp{\na(2\ut+u)}{\na u}_{\Ltom,A}.
\end{align*}
\end{theo}

The functional error estimators $M_+$ and $M_-$ are referred 
as the majorant and the minorant, respectively.
They possess the usual features, e.g., they contain just one constant $\cNa$, which is well known,
and they are sharp. Hence the variational problems for $M_{\pm}$
provide themselves new and equivalent variational formulations
for the system \eqref{diffeq}-\eqref{bceq}. Moreover, since the functions $u$
and vector fields $v$ are at our disposal, one can generate different numerical
schemes to estimate the energy norm of the error $\norm{\na e}_{\Ltom,A}$.
For details see, e.g.,
\cite{MaliRepinNeittaanmaki2013,NeittaanmakiRepin2004,repinbookone} and references therein.

\begin{rem} 
\mylabel{re:majquad}
It is often desirable to have the majorant in the quadratic form
$$M_{+}^2(\na\ut,v)\leq\inf_{\beta>0}
\big(\cNa^2(1+1/\beta)\norm{f+\div v}_{\Ltpoom}^2
+(1+\beta)\norm{\na\ut-A^{-1}v}_{\Ltom,A}^2\big).$$
This form is well suited for computations, since the minimization
with respect to the vector fields $v\in\Dom$ over some finite dimensional subspace
(e.g., generated by finite elements) reduces to solving a system of linear equations.
\end{rem}

\begin{rem} 
\mylabel{re:majind}
If $v\approx A\nabla\uh$, then the first term of the majorant is
close to zero and the second term can be used as an error indicator
to study the distribution of the error over the domain, i.e.,
$$(\na\ut-A^{-1}v)\cdot(A\na\ut-v)\approx\na(\ut-\uh)\cdot A\na(\ut-\uh)\quad\text{ in }\om.$$
The question how to measure the actual performance of the error indicator
(`the accuracy of the symbol $\approx$') is addressed extensively 
in the forthcoming book \cite{MaliRepinNeittaanmaki2013}.
\end{rem}

\section{Non-Conforming A Posteriori Estimates}

To achieve estimates for non-conforming approximations $\ut\notin\Hocmoom+\{u_{0}\}$,
we utilize the simple Helmholtz decomposition
\begin{align}
\mylabel{helm}
\Ltom=\na\Hocmoom\oplus_{A}A^{-1}\Dzom,
\end{align}
where $\Dzom:=\set{v\in\Dom}{\div v=0}$ and $\oplus_{A}$ denotes
the orthogonal sum with respect to the weighted $A$-$\Ltom$-scalar product, see \eqref{defscp}.
The decomposition \eqref{helm} follows immediately from
the projection theorem in Hilbert spaces
$$\Ltom=\ol{\na\Hocom}\oplus_{A}A^{-1}\Dzom$$
and the fact that $\na\Hocmoom=\ol{\na\Hocom}$ is closed in $\Ltom$ by \eqref{poincare}.
Note that the negative divergence
$$-\div:\Dom\subset\Ltom\to\Ltom$$
is the adjoint of the gradient
$$\nac:\Hocom\subset\Ltom\to\Ltom,$$
i.e., $\nac{}^{*}=-\div$, and thus $\Ltom=\ol{R(\nac)}\oplus N(\nac{}^{*})$.
Here, we have used the unweighted standard Sobolev spaces $\Hocom$ and $\Dom$
and the notation $R$ and $N$ for the range and the null space or kernel
of a linear operator, respectively.

Let us now assume that we have an approximation $\vt\in\Ltom$ of $A\na\uh$.
According to \eqref{helm} we decompose the `gradient-error' orthogonally
\begin{align}
\mylabel{Edeco}
\Ltom\ni E:=\na\uh-A^{-1}\vt&=\na\varphi+A^{-1}\psi,\quad\varphi\in\Hocmoom,\,\psi\in\Dzom
\end{align}
and note that it decomposes by Pythagoras' theorem into
\begin{align}
\mylabel{Edecoortho}
\norm{E}_{\Ltom,A}^2&=\norm{\na\varphi}_{\Ltom,A}^2+\norm{A^{-1}\psi}_{\Ltom,A}^2,
\end{align}
which allows us to estimate the two error terms separately.

For $u\in\Hocmoom$ we have by orthogonality
$$\scp{\na\varphi}{\na u}_{\Ltom,A}
=\scp{E}{\na u}_{\Ltom,A}
=\scp{f}{u}_{\Ltom}-\scp{\vt}{\na u}_{\Ltom}.$$
Now we can proceed exactly as in \eqref{calcone} and \eqref{calctwo}
replacing $A\na\ut$ by $\vt$.
More precisely, for all $v\in\Dom$
\begin{align*}
\scp{\na\varphi}{\na u}_{\Ltom,A}
&=\scp{f}{u}_{\Ltom}-\scp{\vt-v+v}{\na u}_{\Ltom}\\
&=\scp{f+\div v}{u}_{\Ltom}-\scp{\vt-v}{\na u}_{\Ltom}
\end{align*}
and hence
\begin{align*}
\scp{\na\varphi}{\na u}_{\Ltom,A}
&\leq\norm{f+\div v}_{\Ltpoom}\norm{u}_{\Ltmoom}
+\norm{A^{-1}(\vt-v)}_{\Ltom,A}\norm{\na u}_{\Ltom,A}\\
&\leq\big(\ub{\cNa\norm{f+\div v}_{\Ltpoom}
+\norm{A^{-1}(\vt-v)}_{\Ltom,A}}_{\ds=M_{+}(A^{-1}\vt,v;f,A)=M_{+}(A^{-1}\vt,v)}\big)
\norm{\na u}_{\Ltom,A}.
\end{align*}
Setting $u:=\varphi$ yields
\begin{align}
\mylabel{nonconfubone}
\norm{\na\varphi}_{\Ltom,A}&\leq\inf_{v\in\Dom}M_{+}(A^{-1}\vt,v).
\end{align}
This estimate is no longer sharp contrary to the conforming case.
We just have the equality $M_{+}(A^{-1}\vt,v)=\norm{E}_{\Ltom,A}$ for $v=A\na\uh$.

For $v\in\Dzom$ and $u\in\Hocmoom+\{u_{0}\}$, we have by orthogonality
and since $\uh-u$ belongs to $\Hocmoom$
\begin{align*}
\scp{A^{-1}\psi}{A^{-1}v}_{\Ltom,A}
&=\scp{E}{A^{-1}v}_{\Ltom,A}\\
&=\ub{\scp{\na(\uh-u)}{v}_{\Ltom}}_{\ds=0}+\scp{\na u-A^{-1}\vt}{A^{-1}v}_{\Ltom,A}\\
&\leq\ub{\norm{\na u-A^{-1}\vt}_{\Ltom,A}}_{\ds=:\tilde{M}_{+}(A^{-1}\vt,\na u;A)=:\tilde{M}_{+}(A^{-1}\vt,\na u)}\norm{A^{-1}v}_{\Ltom,A}.
\end{align*}
Setting $v:=\psi$ yields
\begin{align}
\mylabel{nonconfubtwo}
\norm{A^{-1}\psi}_{\Ltom,A}
&\leq\inf_{u\in\Hocmoom+\{u_{0}\}}\tilde{M}_{+}(A^{-1}\vt,\na u).
\end{align}
Again, this estimate is no longer sharp. We just have
$\tilde{M}_{+}(A^{-1}\vt,\na u)=\norm{E}_{\Ltom,A}$ for $u=\uh$.

For the lower bounds we pick an arbitrary $u\in\Hocmoom$ and compute
\begin{align*}
\norm{\na\varphi}_{\Ltom,A}^2
&\geq2\ub{\scp{\na\varphi}{\na u}_{\Ltom,A}}_{\ds=\scp{E}{\na u}_{\Ltom,A}}
-\norm{\na u}_{\Ltom,A}^2\\
&=2\scp{\na\uh}{\na u}_{\Ltom,A}
-2\scp{A^{-1}\vt}{\na u}_{\Ltom,A}
-\norm{\na u}_{\Ltom,A}^2\\
&=2\scp{f}{u}_{\Ltom}-\scp{\na u+2A^{-1}\vt}{\na u}_{\Ltom,A}\\
&=M_{-}(A^{-1}\vt,u;f,A)=M_{-}(A^{-1}\vt,u).
\end{align*}
Substituting $u=\varphi$ shows that this lower bound is sharp
since $M_{-}(A^{-1}\vt,u)=\norm{\na\varphi}_{\Ltom,A}^2$.

Now we choose $v\in\Dzom$ and $u\in\Hocmoom+\{u_{0}\}$ getting
\begin{align*}
\norm{A^{-1}\psi}_{\Ltom,A}^2
&\geq2\ub{\scp{A^{-1}\psi}{A^{-1}v}_{\Ltom,A}}_{\ds=\scp{E}{A^{-1}v}_{\Ltom,A}}
-\norm{A^{-1}v}_{\Ltom,A}^2\\
&=2\ub{\scp{\na(\uh-u)}{A^{-1}v}_{\Ltom,A}}_{\ds=0}
+2\scp{\na u-A^{-1}\vt}{A^{-1}v}_{\Ltom,A}
-\norm{A^{-1}v}_{\Ltom,A}^2\\
&=\scp{2\na u-A^{-1}(2\vt+v)}{A^{-1}v}_{\Ltom,A}\\
&=:\tilde{M}_{-}(A^{-1}\vt,A^{-1}v,\na u;A)
=:\tilde{M}_{-}(A^{-1}\vt,A^{-1}v,\na u),
\end{align*}
since $\uh-u\in\Hocmoom$.
Also, this second lower bound is still sharp
since we have
$\tilde{M}_{-}(A^{-1}\vt,A^{-1}v,\na u)=\norm{A^{-1}\psi}_{\Ltom,A}^2$
for $v=\psi$ and any $u\in\Hocmoom+\{u_{0}\}$.

\begin{theo}[non-conforming a posteriori estimates]
\mylabel{nonconftheo}
Let $\vt\in\Ltom$. Then
\begin{align*}
&\qquad\norm{\na\uh-A^{-1}\vt}_{\Ltom,A}^2\\
&\leq\calM_{+}(\vt):=\inf_{v\in\Dom}M_{+}^2(A^{-1}\vt,v)
+\inf_{u\in\Hocmoom+\{u_{0}\}}\tilde{M}_{+}^2(A^{-1}\vt,\na u),\\
&\qquad\norm{\na\uh-A^{-1}\vt}_{\Ltom,A}^2\\
&\geq\calM_{-}(\vt):=\sup_{u\in\Hocmoom}M_{-}(A^{-1}\vt,u)
+\sup_{u\in\Hocmoom+\{u_{0}\}}\sup_{v\in\Dzom}
\tilde{M}_{-}(A^{-1}\vt,A^{-1}v,\na u),
\end{align*}
where
\begin{align*}
M_{+}(A^{-1}\vt,v)
&=\cNa\norm{f+\div v}_{\Ltpoom}
+\norm{A^{-1}(\vt-v)}_{\Ltom,A},\\
\tilde{M}_{+}(A^{-1}\vt,\na u)
&=\norm{\na u-A^{-1}\vt}_{\Ltom,A},\\
M_{-}(A^{-1}\vt,u)
&=2\scp{f}{u}_{\Ltom}-\scp{\na u
+2A^{-1}\vt}{\na u}_{\Ltom,A},\\
\tilde{M}_{-}(A^{-1}\vt,A^{-1}v,\na u)
&=\scp{2\na u-A^{-1}(2\vt+v)}{A^{-1}v}_{\Ltom,A}.
\end{align*}
Moreover, as in Remark \ref{re:majquad}
$$M_{+}^2(A^{-1}\vt,v)\leq\inf_{\beta>0}
\big(\cNa^2(1+1/\beta)\norm{f+\div v}_{\Ltpoom}^2
+(1+\beta)\norm{A^{-1}(\vt-v)}_{\Ltom,A}^2\big).$$
\end{theo}

\begin{rem}
\mylabel{nonconfremsharp}
The lower bound is still sharp also in this non-conforming estimate.
As shown before, taking $u=\varphi$ yields
$M_{-}(A^{-1}\vt,u)=\norm{\na\varphi}_{\Ltom,A}^2$
and for $v=\psi$ and arbitrary $u\in\Hocmoom+\{u_{0}\}$ we have
$\tilde{M}_{-}(A^{-1}\vt,A^{-1}v,\na u)=\norm{A^{-1}\psi}_{\Ltom,A}^2$.
Thus, we even have
$$\norm{\na\uh-A^{-1}\vt}_{\Ltom,A}^2
\geq\calM_{-}(\vt)\geq\norm{\na\varphi}_{\Ltom,A}^2
+\norm{A^{-1}\psi}_{\Ltom,A}^2=\norm{\na\uh-A^{-1}\vt}_{\Ltom,A}^2,$$
i.e., $\calM_{-}(\vt)=\norm{\na\uh-A^{-1}\vt}_{\Ltom,A}^2$.
The upper bound might no longer be sharp.
Taking e.g. $u:=\uh$ and $v=A\na\uh$ we get
$$M_{+}(A^{-1}\vt,v)=\tilde{M}_{+}(A^{-1}\vt,\na u)=\norm{\na\uh-A^{-1}\vt}_{\Ltom,A}$$
and thus $\calM_{+}(\vt)\leq2\norm{\na\uh-A^{-1}\vt}_{\Ltom,A}^2$.
So, an overestimation by $2$ is possible.
\end{rem}

\begin{rem}
\mylabel{nonconfremconf}
For conforming approximations $\vt=A\na\ut$, i.e., $A^{-1}\vt=\na\ut$,
with some $\ut\in\Hocmoom+\{u_{0}\}$, we obtain the estimates
from Theorem \ref{conftheo} since
$$\inf_{u\in\Hocmoom+\{u_{0}\}}\tilde{M}_{+}(A^{-1}\vt,\na u)
=\inf_{u\in\Hocmoom+\{u_{0}\}}\norm{\na(u-\ut)}_{\Ltom,A}=0$$
and
$$\sup_{u\in\Hocmoom+\{u_{0}\}}\sup_{v\in\Dzom}
\ub{\tilde{M}_{-}(A^{-1}\vt,A^{-1}v,\na u)}_{\ds=\scp{2\na(u-\ut)-A^{-1}v}{A^{-1}v}_{\Ltom,A}}
=\sup_{v\in\Dzom}-\norm{A^{-1}v}_{\Ltom,A}^2=0$$
because $\scp{\na(u-\ut)}{A^{-1}v}_{\Ltom,A}=\scp{\na(u-\ut)}{v}_{\Ltom}=0$
by $u-\ut\in\Hocmoom$.
\end{rem}

\begin{rem}
\mylabel{nonconfrembderr}
The terms $\tilde{M}_{\pm}$ measure the boundary error.
To see this, let us introduce the scalar trace operator
$\gamma:\Homoom\to\Hgen{1/2}{}{}(\ga)$
and a corresponding extension operator
$\check{\gamma}:\Hgen{1/2}{}{}(\ga)\to\Hoom$.
These are both linear and continuous
(let's say with constants $c_{\gamma}$ and $c_{\check{\gamma}}$)
and $\gamma$ is surjective.
Moreover, $\check{\gamma}$ is a right inverse to $\gamma$.
For an approximation $\vt=A\na\ut$ with $\ut\in\Homoom$ we define
$\check{u}:=\ut+\check{\gamma}\gamma(u_{0}-\ut)$. Then $\gamma\check{u}=\gamma u_{0}$
and hence $\check{u}-u_{0}\in\Hocmoom$. We obtain
\begin{align*}
&\qquad\inf_{u\in\Hocmoom+\{u_{0}\}}
\ub{\tilde{M}_{+}(A^{-1}\vt,\na u)}_{\ds=\norm{\na(u-\ut)}_{\Ltom,A}}\\
&\leq\norm{\na(\check{u}-\ut)}_{\Ltom,A}
=\norm{\na\check{\gamma}\gamma(u_{0}-\ut)}_{\Ltom,A}
\leq c_{\check{\gamma}}\norm{\gamma(u_{0}-\ut)}_{\Hgen{1/2}{}{}(\ga)}
\end{align*}
and since $v\in\Dzom$ by partial integration using the normal trace
$\gamma_{\nu}v\in\Hgen{-1/2}{}{}(\ga)$
\begin{align*}
&\qquad\sup_{u\in\Hocmoom+\{u_{0}\}}\sup_{v\in\Dzom}
\ub{\tilde{M}_{-}(A^{-1}\vt,A^{-1}v,\na u)}_{\ds=2\scp{\na(u-\ut)}{v}_{\Ltom}-\norm{A^{-1}v}^2_{\Ltom,A}}\\
&\geq\sup_{v\in\Dzom}\big(2\scp{\na(\check{u}-\ut)}{v}_{\Ltom}
-\norm{A^{-1}v}_{\Ltom,A}^2\big)\\
&=\sup_{v\in\Dzom}\big(2\scp{\gamma(u_{0}-\ut)}{\gamma_{\nu}v}_{\Hgen{1/2}{}{}(\ga),\Hgen{-1/2}{}{}(\ga)}
-\norm{A^{-1}v}_{\Ltom,A}^2\big).
\end{align*}
\end{rem}

\begin{rem}
\mylabel{remgenbc}
The results of this contribution extend easily and in a canonical way to
exterior elliptic boundary value problems
with pure Neumann or mixed boundary conditions, such as
\begin{align*}
-\div A\nabla u&=f&&\text{in }\om,\\
u&=u_0&&\text{on }\ga_{1},\\
\nu\cdot A\nabla u&=v_{0}&&\text{on }\ga_{2},
\end{align*}
where the boundary $\ga$ decomposes into two parts $\ga_{1}$ and $\ga_{2}$.
\end{rem}

\section{Numerical Tests}

Let $B_{R}:=\set{x\in\rN}{|x|<R}$ denote the open ball,
$E_{R}:=\set{x\in\rN}{|x|>R}$ the exterior domain
and $S_{R}:=\set{x\in\rN}{|x|=R}$ the sphere of radius $R$ 
centered at the origin, respectively, 
as well as $\omR:=\om\cap B_{R}$.
Moreover, let $N=3$, thus $\cN=2$.
In our examples we set $A=\id$, $f=0$ and $u_{0}=1$.
Hence we have $\alpha=1$ and $\cNa=2$.

Therefore, we will consider the exterior Dirichlet Laplace problem, i.e., find
$\uh\in\Homoom$ such that
\begin{align}
\Delta\uh&=0&&\text{in }\om,\label{diffeqnum}\\
\uh&=1&&\text{on }\ga.\label{bceqnum}
\end{align}

It is classical that any solution $u\in\Ltmoom$ (Even $u\in\Ltsom$
with some $s\geq-3/2$ is sufficient.) of $\Delta u=0$ in $E_{R}$
can be represented as as a spherical harmonics expansion with only
negative powers, more precisely as a series of spherical harmonics
$Y_{n,m}$ of order $n$ multiplied by proper powers of the radius
$r^{-(n+1)}$, i.e., for $r>R$
\begin{align}
\mylabel{sphharm}
u_{\Phi}(r,\theta,\varphi)
=\sum_{\substack{n\geq0\\-n\leq m\leq n}}
\zeta_{n,m}r^{-n-1}Y_{n,m}(\theta,\varphi),\quad
\zeta_{n,m}\in\rz,
\end{align}
where $u_{\Phi}:=u\circ\Phi$ and $\Phi$ denotes the usual polar coordinates,
see e.g. \cite{couranthilbertbookone,sobolevbook}.

\begin{rem}
\mylabel{remoneoverr}
For $\om=E_{1}$ the unique solution of \eqref{diffeqnum}-\eqref{bceqnum}
is $\uh=1/r$, which is the first term in the expansion \eqref{sphharm}
corresponding to $n=0$.
We note that even $\uh\in\Ltsom$ as well as $|\na\uh|=1/r^2\in\Ltspoom$ hold
for every $s<-1/2$. For any $1<R<R'$ it is also the unique solution
of the exterior Dirichlet Laplace problem
\begin{align*}
\Delta\uh&=0&&\text{in }E_{R},\\
\uh&=1/R&&\text{on }S_{R}
\intertext{with $\uh\in\Homo(E_{R})$ and of the Dirichlet Laplace problem}
\Delta\uh&=0&&\text{in }B_{R'}\cap E_{R},\\
\uh&=1/R&&\text{on }S_{R},\\
\uh&=1/R'&&\text{on }S_{R'}
\end{align*}
with $\uh\in\Ho(B_{R'}\cap E_{R})$. 
\end{rem}

\begin{figure}[tb]
\begin{center}
\begin{tikzpicture}[scale=2.2]
\draw [rounded corners=10pt] (-1.2,-1.2) rectangle (2.3,1.3);
\draw (-1,1.1) node {$\rN$};
\filldraw [gray!20] (0,0) circle (1); 
\draw (0,0) circle (1);
\draw (1.6,0.4) node {$E_R$};
\draw (1,0) node[anchor=west] {$S_R$};
\filldraw [white] (0,0) circle (0.2); 
\draw (0,0) circle (0.2);
\draw (0.2,0) node[anchor=west] {$\ga$};
\draw (0.3,-0.5) node[anchor=west] {$\omR$};
\draw [<->] (0,0) -- +(45:1);     
\draw (0.35,0.3) node[anchor=west] {$R$};
\draw (-0.6,0.05) node {$\rN\setminus\ol{\om}$};
\draw[->] (-0.4,0.2) parabola bend (-0.2,0.3) (-0.05,0.1);
\end{tikzpicture}
\caption{Domains.} 
\label{fig:domains}
\end{center}
\end{figure}

Of course, the system \eqref{diffeqnum}-\eqref{bceqnum}
is equivalent to find $\uh\in\Homoom$ with $\uh|_{\ga}=1$ and
$$\forall\,\varphi\in\Hocmoom\qquad\scp{\na\uh}{\na\varphi}_{\Ltom}=0,$$
or to minimize the energy
$$\calE(u):=\norm{\na u}_{\Ltom}^2$$
over the set $\set{u\in\Homoom}{u|_{\ga}=1}$.
In order to generate an approximate solution $\ut$, we
split $\om$ into an unbounded and a bounded subdomain, namely
$E_{R}$ and $\omR$, where we pick $R>0$ such that 
$\rt\setminus\om\subset B_{R}$.
The domains are depicted in Figure \ref{fig:domains}.
Our approximation method is based on the assumption
that if $R$ is `large enough',
then the first term of the expansion \eqref{sphharm} dominates 
in the unbounded subdomain $E_R$ and hence we simply assume
from our approximation $\ut\in\Homoom$ the asymptotic behavior
$$\ut|_{E_{R}}=\frac{\zeta}{r},$$ 
where $\zeta$ is an unknown real constant.
In the bounded subdomain $\omR$, 
the approximation $\utR:=\ut|_{\omR}\in\HoomR$ 
must satisfy the boundary condition
$\utR|_{\ga}=1$ and the continuity condition $\utR|_{S_{R}}=\zeta/R$
to ensure $\ut\in\Homoom$.
Then, for approximations of the prescribed type, the problem 
is reduced to minimize the energy
$$\calE(\utRz,\zeta):=\calE(\ut)=\norm{\na\ut}_{\Ltom}^2
=\norm{\na\utRz}_{\Lt(\omR)}^2+\frac{4\pi\zeta^2}{R}$$
with respect to $\zeta\in\rz$ and $\utRz\in\HoomR$ 
with $\utRz|_{\ga}=1$ and $\utRz|_{S_{R}}=\zeta/R$.
We propose an iteration procedure to minimize the quadratic energy or functional $\calE$, 
which is described in Algorithm \ref{alg1}. It is based on the decomposition 
$$\utRz=\utRzz+\utRo+\frac{\zeta}{R}\utRt,$$
where $\utRzz\in\HocomR$ and $\utRo,\utRt\in\HoomR$ satisfy certain boundary conditions, i.e.,
\begin{align}
\mylabel{utRbc}
\utRo|_{\ga}=1,\quad\utRo|_{S_{R}}=0,\qquad\utRt|_{\ga}=0,\quad\utRt|_{S_{R}}=1,
\end{align}
as well as
\begin{align}
\mylabel{utRLaplacenotzero}
\Delta\utRt\neq0.
\end{align}
The two functions $\utRo,\utRt$ take only care of the boundary conditions
and are fixed during the iteration procedure.
This means that the energy $\calE(\utRz,\zeta)$ is minimized with respect to 
$\zeta\in\rz$ and $\utRzz$ in $\HocomR$.
We note that \eqref{utRLaplacenotzero} is crucial for the iteration process
since otherwise the update $\ut_{k}$ in Algorithm \ref{alg1}
would not depend on the previous $\zeta_{k}$.

\begin{algorithm}[h!]
\caption{\sf\quad Minimization of the energy $\calE(\utRz,\zeta)$} 
\mylabel{alg1}
\begin{itemize}
\item{\sf step 0:} 
Pick any $\utRo,\utRt\in\HoomR$ with \eqref{utRbc} and \eqref{utRLaplacenotzero}
and set $k:=1$ and $\zeta_{k}:=1$.
\item{\sf step 1:} 
Minimize the quadratic energy 
$$\calE_{k}(u):=\calE(u+\utRo+\frac{\zeta_{k}}{R}\utRt,\zeta_{k})$$
with respect to $u\in\HocomR$, i.e., find $u\in\HocomR$ such that
$$\forall\,\varphi\in\HocomR\quad
\scp{\na u}{\na\varphi}_{\Lt(\omR)}
=-\scp{\na\utRo}{\na\varphi}_{\Lt(\omR)}
-\frac{\zeta_{k}}{R}\scp{\na\utRt}{\na\varphi}_{\Lt(\omR)}.$$
Set $\utk:=u$.
\item{\sf step 2:} 
Minimize the second order polynomial 
$$p_{k}(\zeta):=\calE(\utk+\utRo+\frac{\zeta}{R}\utRt,\zeta)
=\norm{\na\utk+\na\utRo+\frac{\zeta}{R}\na\utRt}_{\Lt(\omR)}^2
+\frac{4\pi\zeta^2}{R}$$
with respect to $\zeta$, i.e., set
$\ds\zeta:=-R\frac{\scp{\na(\utk+\utRo)}{\na\utRt}_{\Lt(\omR)}}{4\pi R+\norm{\na\utRt}_{\Lt(\omR)}^2}$.\\
Set $\zeta_{k+1}:=\zeta$.
\item{\sf step 3:} 
Set $k:=k+1$ and return to step 1, 
unless $|\zeta_{k}-\zeta_{k-1}|/|\zeta_{k}|$ is small.
\item{\sf step 4:} 
Set $\utRzz:=\ut_{k-1}$, $\zeta:=\zeta_{k}$ and
$$\ut:=\begin{cases}\ds\utRz:=\utRzz+\utRo+\frac{\zeta}{R}\utRt&\text{ in }\omR\\
\ds\frac{\zeta}{r}&\text{ in }E_{R}\end{cases}.$$
\end{itemize}
\end{algorithm}

The conforming estimates from Theorem \ref{conftheo} 
involving the free variables $u\in\Hocmoom$ and $v\in\Dom$
are used to estimate the approximation error. 
For the variable $u$ we simply choose $u_{R}\in\HocomR$ 
and extend $u_{R}$ by zero to $\om$, which defines a proper $u$.
To restrict all computations to $\omR$, the best choice for $v$ in $E_{R}$ is
$v|_{E_{R}}:=\na\ut|_{E_{R}}=-\zeta r^{-2}e^r$ 
with the unit radial vector $e^r(x):=x/|x|$.
Picking vector fields $v_{R}$ as restrictions from $\Dom$ to $\omR$, i.e., 
$v_{R}\in\DomR$, we need that the extensions
\begin{align}
\mylabel{defuv}
v:=\begin{cases}v_{R}&\text{in }\omR\\-\zeta r^{-2}e^r&\text{in }E_{R}\end{cases}
\end{align}
belong to $\Dom$. Hence, $e^r\cdot v$, the normal component of $v$,
must be continuous across $S_{R}$. 
Therefore, we get on $S_{R}$ the transmission condition $e^R\cdot v_{R}=-\zeta/R^2$.
Thus, any $v_{R}$ in
$$\VomR:=\set{\psi\in\DomR}{e^R\cdot\psi|_{S_{R}}=-\zeta/R^2},$$
which is extended by \eqref{defuv} to $\om$ belongs to $\Dom$.
We note that then $\na\ut=v$ and $\div v=\Delta\ut=0$ holds in $E_{R}$.
Now, the estimates from Theorem \ref{conftheo} read
\begin{align*}
&\qquad\norm{\na e}_{\Ltom}^2=\norm{\na(\uh-\ut)}_{\Ltom}^2\\
&\leq\inf_{v\in\Dom}\inf_{\beta>0}
\big(4(1+1/\beta)\norm{\div v}_{\Ltpoom}^2
+(1+\beta)\norm{\na\ut-v}_{\Ltom}^2\big)\\
&=\inf_{v\in\Dom}\inf_{\beta>0}
\big(4(1+1/\beta)\int_{\om}(1+r^2)|\div v|^2\dl
+(1+\beta)\int_{\om}|\na\ut-v|^2\dl\big)\\
&\leq\inf_{v_{R}\in\VomR}\inf_{\beta>0}
\ub{\big(4(1+\beta)\int_{\omR}(1+r^2)|\div v_{R}|^2\dl
+(1+\beta)\int_{\omR}|\na\ut-v_{R}|^2\dl\big)}_{\ds=:M_{+,R,\beta}^2(\na\ut,v_{R})}
\intertext{and}
&\qquad\norm{\na e}_{\Ltom}^2=\norm{\na(\uh-\ut)}_{\Ltom}^2\\
&\geq\sup_{u\in\Hocmoom}
-\scp{\na(2\ut+u)}{\na u}_{\Ltom}
=\sup_{u\in\Hocmoom}
-\int_{\om}\na(2\ut+u)\cdot\na u\dl\\
&\geq\sup_{u_{R}\in\HocomR}
\ub{-\int_{\omR}\na(2\ut+u_{R})\cdot\na u_{R}\dl}_{\ds=:M_{-,R}(\na\ut,u_{R})}.
\end{align*}
Therefore, we have reduced the computations of the lower and upper bounds
to minimization problems taking place only in the bounded domain $\omR$.

\begin{ex}
\mylabel{ex1}
We set $\om:=E_{1}$, i.e., $\om$ is the exterior of the closed unit ball.  
As already mentioned,
this problem admits the unique solution $\uh=1/r$.
Hence, we know $\zeta=1$ a priori, but ignore it for the computations.
We use the symmetry of the problem
and thus our computations take place in just one octant of $\omR$.
The mesh and the computational domain are depicted in Figure \ref{oct1}.
The resulting relative errors and error estimates
are presented in Table \ref{tab1} for three different values of $R$.
Additionally, we study the error indicator generated by the majorant, see Remark \ref{re:majind}.
The exact error contribution over an element $T$ is
$\norm{\na e}_{\Lt(T)}^2$
and the one indicated by the second term of the majorant is
$\norm{\na\ut-v_{R}}_{\Lt(T)}^2$, where $v_{R}$ is obtained 
via the minimization of the majorant. Both quantities
are depicted on the plane $x_{1}+x_{3}=0$ in Figure \ref{fig:errind1}.
\begin{figure}[tb]
\centering
\includegraphics[width=7.5cm,clip=true,trim=3cm 2cm 3cm 2cm]{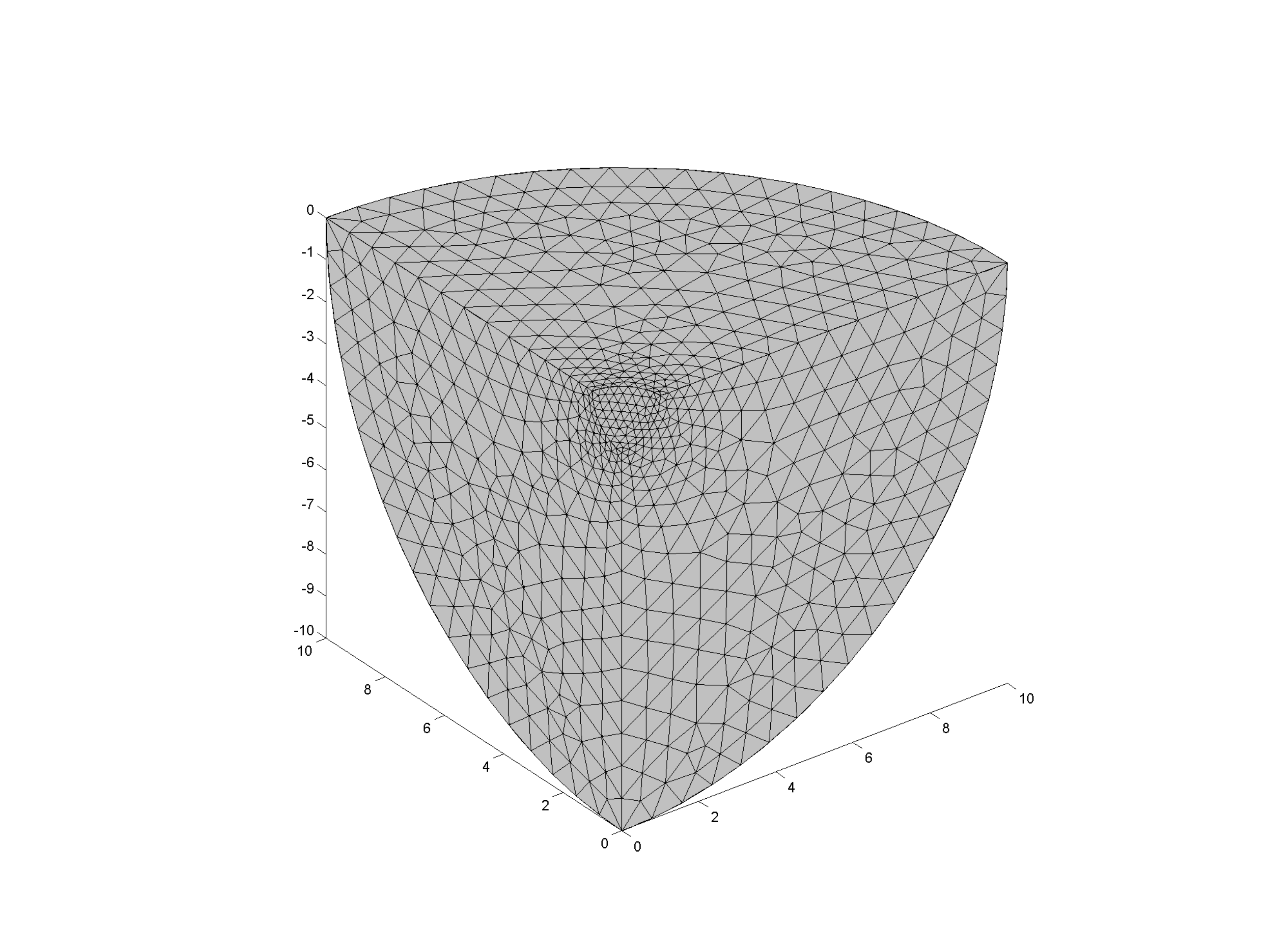}
\includegraphics[width=7.5cm,clip=true,trim=3cm 2cm 3cm 2cm]{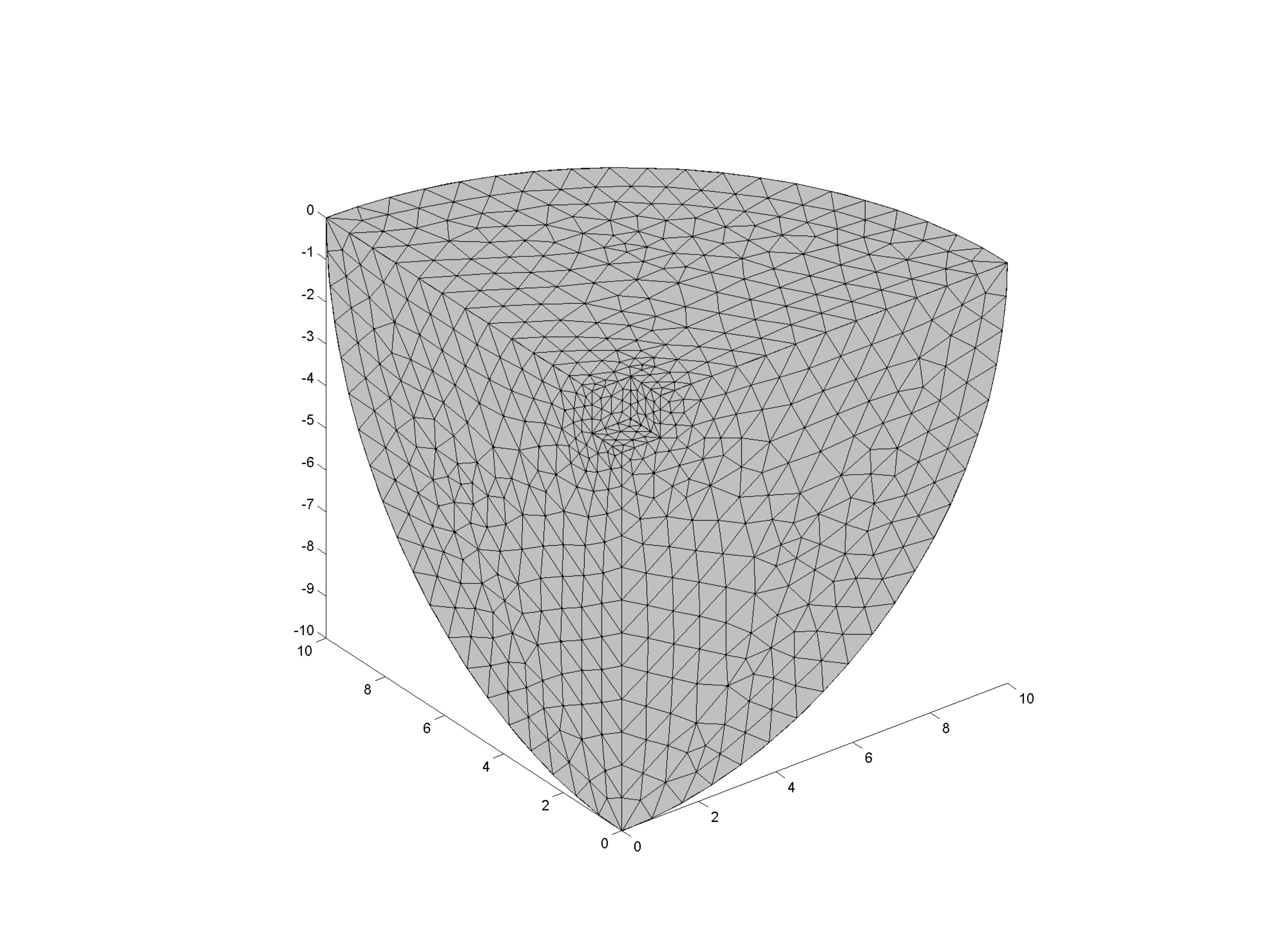}
\caption{Computational domains and meshes (one octant of $\omR$, where $R=10$)
for Example \ref{ex1} (left) and Example \ref{ex2} (right).}
\mylabel{oct1}
\end{figure}
\begin{table}
\centering
\caption{Example \ref{ex1}, the exact error and the computed error estimates.}
\mylabel{tab1}
\begin{tabular}{@{}rrccc@{}}
\hline
$R$ & 
$N_t$ & 
$\frac{M_{+,R,\beta}^2(\na\ut,v_{R})}{\norm{\na\ut}_{\Ltom}^2}$ &
$\frac{\norm{\na e}_{\Ltom}^2}{\norm{\na\ut}_{\Ltom}^2}$ &
$\frac{M_{-,R}(\na\ut,u_{R})}{\norm{\na\ut}_{\Ltom}^2}$ \\
\hline\hline
 5 &  1926 & 4.57\% & 3.75\% & 3.61\% \\
 5 &  6223 & 1.39\% & 1.31\% & 1.29\% \\
 5 & 17742 & 0.63\% & 0.62\% & 0.61\% \\
\hline
10 &  2214 & 7.19\% & 5.68\% & 5.43\% \\
10 &  5264 & 2.82\% & 2.55\% & 2.51\% \\
10 & 11506 & 1.80\% & 1.68\% & 1.65\% \\
\hline
20 &  3446 & 7.26\% & 5.73\% & 5.47\% \\
20 &  6676 & 2.95\% & 2.64\% & 2.60\% \\
20 & 13370 & 1.86\% & 1.74\% & 1.71\% \\
\hline
\end{tabular}
\end{table}
\begin{figure}
\centering
\includegraphics[width=7.5cm,clip=true,trim=3cm 2cm 3cm 2cm]{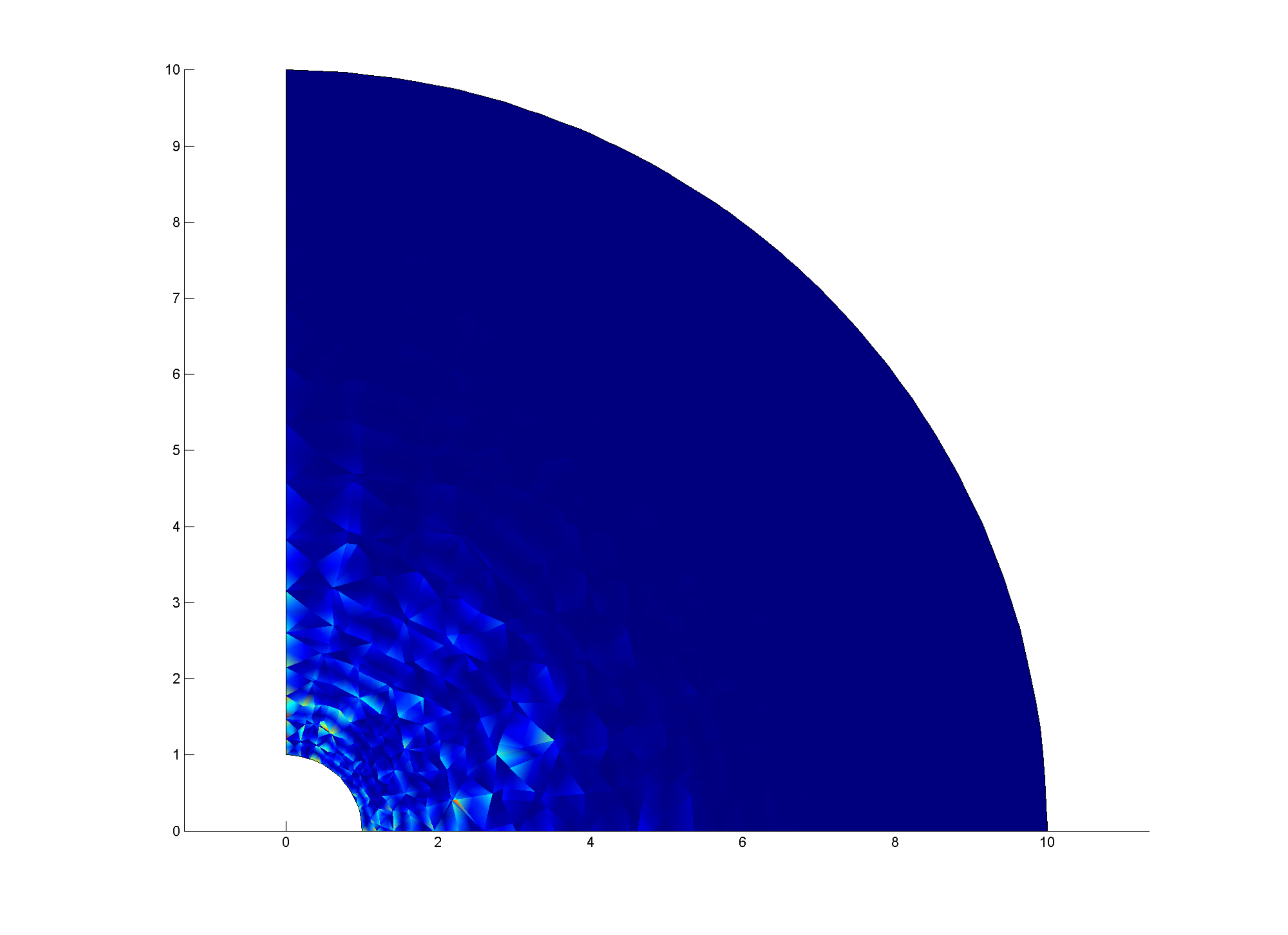}
\includegraphics[width=7.5cm,clip=true,trim=3cm 2cm 3cm 2cm]{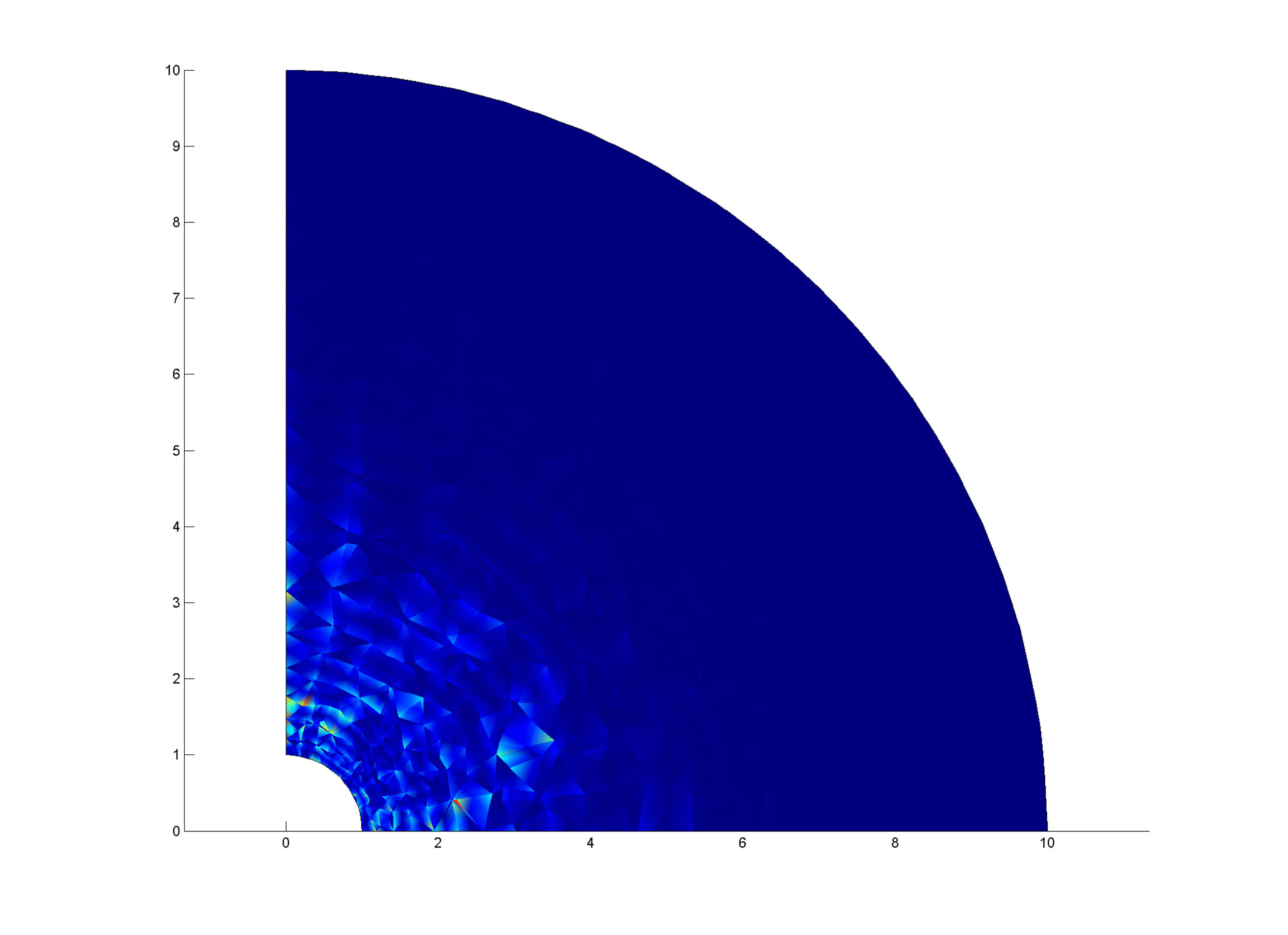}
\caption{Example \ref{ex1}, true error $\norm{\na e}_{\Lt(T)}^2$ (left)
and error indicator $\norm{\na\ut-v_{R}}_{\Lt(T)}^2$ (right).}
\mylabel{fig:errind1}
\end{figure}
\end{ex}

The boundary value problem in step 1 of Algorithm \ref{alg1} 
was solved by the finite element method.
We applied first order nodal tetrahedral elements, 
where the mesh was constructed by Comsol 4.3
and the emerging system of linear equations 
was solved using a standard Matlab solver.
When minimizing the majorant $M_{+,R,\beta}^2(\na\ut,v_{R})$ 
with respect to $v_{R}\in\VomR$, we applied second order tetrahedral 
finite elements for each component of $v_{R}$ 
with $v_{R}|_{S_{R}}=-\zeta R^{-2}e^R$ and hence $e^R\cdot v_{R}|_{S_{R}}=-\zeta/R^{2}$.
Similarly, the minorant $M_{-,R}(\na\ut,u_{R})$ 
was maximized using second order tetrahedral finite elements for $u_{R}\in\HocomR$.
A natural choice is to use $\Sobolev(\div)$-conforming 
Raviart-Thomas elements \cite{RaviartThomas1977} to compute $v_{R}$.
However, $\Ho$-elements are applicable for smooth problems.

\begin{ex}
\mylabel{ex2}
We set $\om:=\rt\setminus[-1,1]^3$, i.e., $\om$ is the exterior of a closed cube. 
For this problem, the exact solution $\uh$ is not known. 
The octant of $\omR$ for $R=10$ used in computations is depicted in Figure \ref{oct1}.
In Algorithm \ref{alg1}, we selected $\utRo$ and $\utRt$ 
as finite element solutions of Dirichlet Laplace problems
with proper boundary conditions.
The respective error bounds are presented in Table \ref{tab2}.
\begin{table}
\centering
\caption{Error bounds for Example \ref{ex2}.}
\begin{tabular}{@{}rcc@{}}
\hline
$N_t$ & 
$\frac{M_{+,R,\beta}^2(\na\ut,v_{R})}{\norm{\na\ut}_{\Ltom}^2}$ &
$\frac{M_{-,R}(\na\ut,u_{R})}{\norm{\na\ut}_{\Ltom}^2}$ \\
\hline\hline
 2513 & 12.52\% & 4.65\% \\
 5015 & 9.48\%  & 2.91\% \\
13772 & 7.08\%  & 1.63\% \\
\hline
\end{tabular}
\mylabel{tab2}
\end{table}
\end{ex}

These examples show that functional a posteriori error estimates 
provide two-sided bounds of the error.
Of course, the accuracy depends on the method
used to generate the free variables $v$ and $u$ 
in the majorant and minorant, respectively. 
The applied methods should be selected balancing 
the desired accuracy of the error estimate 
and the computational expenditures.

\begin{acknow}
We heartily thank Sergey Repin for his continuous support
and the many nice and enlightening discussions.
\end{acknow}

\bibliographystyle{plain}
\bibliography{/Users/paule/Library/texmf/tex/TeXinput/bibtex/paule}

\appendix
\section{Appendix: Easier but Weaker Estimates}

We want to point out that we can prove
a variant of Theorem \ref{nonconftheo}
by another, much simpler technique using just the triangle inequality
instead of the Helmholtz decomposition.
The drawbacks are that we get for the upper bound a factor larger than $1$, 
e.g. $5$, which overestimates a bit more,
and for the lower bound we miss one term.
To see this, let $u\in\Hocmoom+\{u_{0}\}$. Then
$$\norm{E}_{\Ltom,A}
=\norm{\na\uh-A^{-1}\vt}_{\Ltom,A}
\leq\norm{\na(\uh-u)}_{\Ltom,A}
+\norm{\na u-A^{-1}\vt}_{\Ltom,A},$$
and we can further estimate by Theorem \ref{conftheo} with $\ut=u$
for any $v\in\Dom$
\begin{align*}
\norm{E}_{\Ltom,A}
&\leq\ub{\cNa\norm{f+\div v}_{\Ltpoom}
+\norm{\na u-A^{-1}v}_{\Ltom,A}}_{\ds=M_{+}(\na u,v)}
+\norm{\na u-A^{-1}\vt}_{\Ltom,A}\\
&\leq\ub{\cNa\norm{f+\div v}_{\Ltpoom}
+\norm{A^{-1}(\vt-v)}_{\Ltom,A}}_{\ds=M_{+}(A^{-1}\vt,v)}
+2\ub{\norm{\na u-A^{-1}\vt}_{\Ltom,A}}_{\ds=\tilde{M}_{+}(A^{-1}\vt,\na u)}.
\end{align*}
Therefore, we get the same upper bound but with less good factors, i.e.,
for any $\theta>0$
\begin{align*}
\norm{E}_{\Ltom,A}^2
&\leq\tilde{\calM}_{+}(\vt):=(1+4/\theta)\inf_{v\in\Dom}M_{+}^2(A^{-1}\vt,v)\\
&\qquad\qquad\qquad\qquad+(4+\theta)\inf_{u\in\Hocmoom+\{u_{0}\}}\tilde{M}_{+}^2(A^{-1}\vt,\na u)
\end{align*}
with $\tilde{\calM}_{+}(\vt)\geq\calM_{+}(\vt)$.
E.g. for $\theta=1$ we have $1+4/\theta=4+\theta=5$
and $\tilde{\calM}_{+}(\vt)=5\calM_{+}(\vt)$.
For the lower bound and $u\in\Hocmoom$ we simply have
\begin{align*}
\norm{E}_{\Ltom,A}^2
&=\norm{\na\uh-A^{-1}\vt}_{\Ltom,A}^2\\
&\geq2\scp{\na\uh-A^{-1}\vt}{\na u}_{\Ltom,A}-\norm{\na u}_{\Ltom,A}^2\\
&=2\scp{f}{u}_{\Ltom}-\scp{\na u+2A^{-1}\vt}{\na u}_{\Ltom,A}
=M_{-}(A^{-1}\vt,u)
\end{align*}
and the additional term $\tilde{M}_{-}(A^{-1}\vt,A^{-1}v,\na u)$ 
does not appear. Thus,
$$\norm{E}_{\Ltom,A}^2
\geq\tilde{\calM}_{-}(\vt):=\sup_{u\in\Hocmoom}M_{-}(A^{-1}\vt,u)$$
with $\tilde{\calM}_{-}(\vt)\leq\calM_{-}(\vt)$.

\end{document}

%% file: head_paper.tex
\usepackage{a4,exscale,ifthen,amsfonts,amssymb,amsmath,amscd,graphicx}
\usepackage[english]{babel}
\usepackage[mathscr]{eucal}

\author{{\sf\pauthor}}
\numberwithin{equation}{section}

\newenvironment{acknow}{{\vspace*{1cm}\noindent\bf Acknowledgements }}{}
\newcommand{\bewboxw}{\mbox{}\hfill $\square$ \\}

\newcommand{\keywords}[1]{{\noindent\bf Key Words }#1}

\ifthenelse{\equal{\mylabelonoff}{on}}
{\newcommand{\mylabel}[1]{\label{#1}\fbox{{\rm #1}}}}{\newcommand{\mylabel}[1]{\label{#1}\makebox[0mm][]{}}}
\ifthenelse{\equal{\allowdisbrk}{yes}}
{\allowdisplaybreaks}{}

%% file: head_paule_olli.tex
\usepackage[mathscr]{eucal}
\usepackage{color}
\usepackage{tikz}
\usepackage{fancyhdr}
\usepackage{caption}
\usepackage{array,multirow}
\usepackage{multicol}
\usepackage{booktabs}
\usepackage{algorithm}
\usepackage{algorithmic}
\usepackage[all]{xy}

\newcommand{\ds}{\displaystyle}
\newcommand{\ol}{\overline}

\newcommand{\ub}{\underbrace}

\newcommand{\rz}{\mathbb{R}}

\newcommand{\rt}{\rz^3}

\newcommand{\rN}{\rz^N}
\newcommand{\rNtN}{\rz^{N\times N}}

\DeclareMathOperator{\p}{\partial}

\newcommand{\na}{\nabla}

\newcommand{\nac}{\overset{\circ}{\na}}

\renewcommand{\div}{\operatorname{div}}

\newcommand{\om}{\Omega}
\newcommand{\pom}{\p\!\om}

\newcommand{\ga}{\Gamma}

\DeclareMathOperator{\id}{id}
\newcommand{\dl}{\,d\lambda}

\def\set#1#2{\{#1\,:\,#2\}}

\newcommand{\cN}{c_{N}}
\newcommand{\cNa}{c_{N,\alpha}}

\def\ut{\tilde{u}}
\def\uh{\hat{u}}
\def\vt{\tilde{v}}
\def\utk{\tilde{u}_{k}}
\def\utR{\tilde{u}_{R}}
\def\utRo{\tilde{u}_{R,1}}
\def\utRt{\tilde{u}_{R,2}}

\def\utRz{\tilde{u}_{R,\zeta}}
\def\utRzz{\tilde{u}_{R,\zeta,0}}

\def\calM{\mathcal{M}}
\def\calE{\mathcal{E}}

\DeclareMathOperator{\Lebesgue}{\mathsf{L}}
\newcommand{\Lgen}[2]{\Lebesgue^{#1}_{#2}}
\def\Li{\Lgen{\infty}{}}
\def\Liom{\Li(\om)}
\def\Lo{\Lgen{1}{}}
\def\Loom{\Lo(\om)}
\def\Lt{\Lgen{2}{}}
\def\Ltom{\Lt(\om)}

\def\Lts{\Lgen{2}{s}}
\def\Ltspo{\Lgen{2}{s+1}}

\def\Ltmo{\Lgen{2}{-1}}
\def\Ltz{\Lgen{2}{0}}
\def\Ltpo{\Lgen{2}{1}}

\def\Ltsom{\Lts(\om)}
\def\Ltspoom{\Ltspo(\om)}

\def\Ltmoom{\Ltmo(\om)}
\def\Ltzom{\Ltz(\om)}
\def\Ltpoom{\Ltpo(\om)}

\DeclareMathOperator{\Sobolev}{\mathsf{H}}
\newcommand{\Hgen}[3]{\overset{#3}{\Sobolev}{}^{#1}_{#2}}
\def\Ho{\Hgen{1}{}{}}

\def\Hoom{\Ho(\om)}

\def\Hoc{\Hgen{1}{}{\circ}}
\def\Hocom{\Hoc(\om)}

\def\Homo{\Hgen{1}{-1}{}}
\def\Hocmo{\Hgen{1}{-1}{\circ}}
\def\Homoom{\Homo(\om)}
\def\Hocmoom{\Hocmo(\om)}

\DeclareMathOperator{\Cont}{\mathsf{C}}
\newcommand{\Cgen}[2]{\overset{#2}{\Cont}{}^{#1}}

\def\Cic{\Cgen{\infty}{\circ}}

\def\Cicom{\Cic(\om)}

\newcommand{\scp}[2]{\left\langle#1,#2\right\rangle}

\newcommand{\norm}[1]{\left|\normdst\left|#1\right|\normdst\right|}

\newtheorem{lem}{Lemma}

\newtheorem{theo}[lem]{Theorem}

\newtheorem{rem}[lem]{Remark}

\newtheorem{ex}{Example}

\DeclareMathOperator{\divspace}{\mathsf{D}}
\newcommand{\divgen}[2]{\overset{#2}{\divspace}{}_{#1}}
\newcommand{\divgenom}[2]{\divgen{#1}{#2}(\om)}
\newcommand{\dom}{\divgenom{}{}}

\newcommand{\dzom}{\divgenom{0}{}}

\def\Dom{\dom}
\def\Dzom{\dzom}
\def\omR{\om_{R}}
\def\DomR{\divgen{}{}(\omR)}

\def\HoomR{\Ho(\omR)}
\def\HocomR{\Hoc(\omR)}

\DeclareMathOperator{\V}{\mathsf{V}}

\def\VomR{\V(\omR)}

\renewcommand{\norm}[1]{\left|#1\right|}